\documentclass{article}

\usepackage{xypic}
\usepackage{a4}
\usepackage{amsmath,amssymb}

\newtheorem{theorem}{Theorem}[section]
\newtheorem{proposition}[theorem]{Proposition}

\newtheorem{definition}[theorem]{Definition}
\newenvironment{proof}{\smallskip\par{\sc Proof.}\enspace}%
{{\unskip\nobreak\hfil\penalty50\hskip2em                 %
         \hbox{}\nobreak\hfil{\rule[-1pt]{4pt}{8pt}}      %
         \parfillskip=0pt\finalhyphendemerits=0           %
         \par\medskip}}                                   %

\DeclareMathAlphabet{\eufrak}{U}{}{}{}  
\SetMathAlphabet\eufrak{normal}{U}{euf}{m}{n}
\SetMathAlphabet\eufrak{bold}{U}{euf}{b}{n}

\begin{document}

\title{What is Stochastic Independence?}

\author{Uwe Franz}

\date{
Institut f\"ur Mathematik und Informatik \\
Ernst-Moritz-Arndt-Universit\"at Greifswald \\
Friedrich-Ludwig-Jahn-Str.\ 15 a \\
 D-17487 Greifswald, Germany \\
Email: franz@uni-greifswald.de \\
http://hyperwave.math-inf.uni-greifswald.de/algebra/franz}
\maketitle

\abstract{
The notion of a tensor product with projections or with inclusions is defined. It is shown that the definition of stochastic independence relies on such a structure and that independence can be defined in an arbitrary category with a tensor product with inclusions or projections. In this context, the classifications of quantum stochastic independence by Muraki, Ben Ghorbal, and Sch\"urmann become classifications of the tensor products with inclusions for the categories of algebraic probability spaces and non-unital algebraic probability spaces. The notion of a reduction of one independence to another is also introduced.  As examples the reductions of Fermi independence and boolean, monotone, and anti-monotone independence to tensor independence are presented.
}

\section{Introduction}

In this paper we will deal with the question, what stochastic independence is. Since the work of Speicher\cite{speicher96} and Sch\"urmann\cite{schuermann95,benghorbal01,benghorbal+schuermann99} we know that a `universal' notion of independence should come with a product that allows to construct the joint distribution of two independent random variables from their marginal distributions. It turned out that in classical probability there exists only one such product satisfying a natural set of axioms. But there are several different good notions of independence in non-commutative probability. The most important ones were classified in the work of Speicher\cite{speicher96}, Ben Ghorbal and Sch\"urmann\cite{benghorbal+schuermann99,benghorbal01}, and Muraki\cite{muraki01b,muraki02}, they are tensor independence, free independence, boolean independence, monotone independence and anti-monotone
independence.

We present and motivate here the axiomatic framework used in these articles. We show that the classical notion of stochastic independence is based on a kind of product in the category of probability spaces, which is intermediate to the notion of a (universal) product in category theory - which does not exist in this category - and the notion of a tensor product. Furthermore, we show that the classification of stochastic independence by Ben Ghorbal and Sch\"urmann\cite{schuermann95,benghorbal+schuermann99,benghorbal01} and by Muraki\cite{muraki01b,muraki02} is also based on such a product, which we call a tensor product with projections or inclusions, cf.\ Definition \ref{def-tensor-cat-with}. This notion allows to define independence for arbitrary categories, see Definition \ref{def-independence}. If independence is something that depends on a tensor product and projections or inclusions between the original objects and their tensor product, then it is clear that a map between categories that preserves independence should be tensor functor with an additional structure that takes care of the projections or inclusions. This is formalized in Definition \ref{def-reduction}. We show in several examples that these notions are really the correct one, see Subsection \ref{example-class-indep}, \ref{tensor indep}, \ref{fermi-indep}, and \ref{reduction-fermi-bose}, and Section \ref{section-reduction-boolean-etc}.

But let us first look at the notion of independence in classical probability.

\section{Independence for Classical Random Variables}

Two random variables $X_1:(\Omega,\mathcal{F},P)\to(E_1,\mathcal{E}_1)$ and
$X_2:(\Omega,\mathcal{F},P)\to(E_2,\mathcal{E}_2)$, defined on the same
probability space $(\Omega,\mathcal{F},P)$ and with values in two possibly
distinct measurable spaces $(E_1,\mathcal{E}_1)$ and $(E_2,\mathcal{E}_2)$,
are called {\em stochastically independent}\index{independent!stochastically} (or simply {\em independent}) w.r.t.\ $P$, if the $\sigma$-algebras $X_1^{-1}(\mathcal{E}_1)$ and $X_2^{-1}(\mathcal{E}_2)$ are independent w.r.t.\ $P$, i.e.\ if
\[
P\big((X^{-1}_1(M_1)\cap X^{-1}_2(M_2)\big) =P\big((X^{-1}_1(M_1)\big)P\big( X^{-1}_2(M_2)\big) 
\]
holds for all $M_1\in\mathcal{E}_1$, $M_2\in\mathcal{E}_2$. If there is no
danger of confusion, then the reference to the measure $P$ is often omitted.

This definition can easily be extended to arbitrary families of random
variables. A family
$\big(X_j:(\Omega,\mathcal{F},P)\to(E_j,\mathcal{E}_j))_{j\in J}$, indexed by
some set $J$, is called independent, if
\[
P\left(\bigcap_{k=1}^n(X^{-1}_{j_k}(M_{j_k})\right) = \prod_{k=1}^n
P\big(X^{-1}_{j_k}(M_{j_k})\big)
\]
holds for all $n\in\mathbb{N}$ and all choices of indices $k_1,\ldots,k_n\in
J$ with $j_k\not=j_\ell$ for $j\not=\ell$, and all choices of measurable sets
$M_{j_k}\in\mathcal{E}_{j_k}$.

There are many equivalent formulations for independence, consider, e.g., the
following proposition.

\begin{proposition}
Let $X_1$ and $X_2$ be two real-valued random variables. The following are
equivalent.
\begin{itemize}
\item[(i)]
$X_1$ and $X_2$ are independent.
\item[(ii)] For all bounded measurable functions $f_1,f_2$ on $\mathbb{R}$ we
  have
\[
\mathbb{E}\big(f_1(X_1)f_2(X_2)\big)=
\mathbb{E}\big(f_1(X_1)\big)\mathbb{E}\big(f_2(X_2)\big).
\]
\item[(iii)]
The probability space $(\mathbb{R}^2,\mathcal{B}(\mathbb{R}^2),P_{(X_1,X_2)})$ is the
product of the probability spaces
$(\mathbb{R},\mathcal{B}(\mathbb{R}),P_{X_1})$ and
$(\mathbb{R},\mathcal{B}(\mathbb{R}),P_{X_2})$, i.e.
\[
P_{(X_1,X_2)} = P_{X_1}\otimes P_{X_2}.
\]
\end{itemize}
\end{proposition}

We see that stochastic independence can be reinterpreted as a rule to compute
the joint distribution of two random variables from their marginal
distribution. More precisely, their joint distribution can be computed as a
product of their marginal distributions. This product is associative and can
also be iterated to compute the joint distribution of more than two independent
random variables.

The classifications of independence for non-commutative probability
\cite{speicher96,benghorbal+schuermann99,benghorbal01,muraki01b,muraki02}
that we are interested in are based on redefining independence as a
product satisfying certain natural axioms.

\section{Tensor Categories and Independence}

We will now define the notion of independence in the language of category theory. The usual notion of independence for classical probability theory and the independences classified in \cite{speicher96,benghorbal+schuermann99,benghorbal01,muraki01b,muraki02} will then be instances of this general notion obtained by considering the category of classical probability spaces or the category of algebraic probability spaces.

First we recall the definitions of a product, coproduct and a tensor product, see also MacLane\cite{maclane98} for a more detailed introduction. Then we introduce tensor categories with inclusions or projections. This notion is weaker than that of a product or coproduct, but stronger than that of a tensor category. It is exactly what we need to get an interesting notion of independence.

\begin{definition} \label{def-product} (See, e.g., Maclane\cite{maclane98})
A tuple $(B_1\,\Pi\, B_2,\pi_1,\pi_2)$ is called a {\em product} or {\em universal product} of the objects $B_1$ and $B_2$ in the category $\mathcal{C}$, if for any object $A\in{\rm Ob}\,\mathcal{C}$ and any morphisms $f_1:A\to B_1$ and $f_2:A\to B_2$ there exists a unique morphism $h$ such that the following diagram commutes,
\[
\xymatrix{
& A \ar[dl]_{f_1} \ar@{.>}[d]|-h \ar[dr]^{f_2} & \\
B_1 & B_1\,\Pi\, B_2 \ar[l]^{\pi_1} \ar[r]_{\pi_2} & B_2.
}
\]
An object $K$ is called {\em terminal}, if for all objects $A\in{\rm Ob}\,\mathcal{C}$ there exists exactly one morphism from $A$ to $K$.
\end{definition}
The product of two objects is unique up to isomorphism, if it exists. Furthermore, the operation of taking products is commutative and associative up to isomorphism and therefore, if a category has a terminal object and a product for any two objects, then one can also define a product for any finite set of objects.

The notion of coproduct is dual to that of a product, i.e., its defining property can be obtained from that of the product by `reverting the arrows'. The notion dual to terminal object is an initial object, i.e.\ an object $K$ such that for any object $A$ of $\mathcal{C}$ there exists a unique morphism from $K$ to $A$.

Let us now recall the definition of a tensor category.
\begin{definition}
A category $(\mathcal{C},\square)$ equipped with a bifunctor $\square:\mathcal{C}\times\mathcal{C}\to\mathcal{C}$, called {\em tensor product}, that is associative up to a natural isomorphism
\[
\alpha_{A,B,C}:A\square(B\square C) \stackrel{\cong}{\to} (A\square B)\square C,\qquad \text{ for all } A,B,C\in{\rm Ob}\,\mathcal{C},
\]
 and an element $E$ that is, up to isomorphisms
\[
\lambda_A:E\square A\stackrel{\cong}{\to}A, \quad\text{ and }\quad \rho_A:A\square E\stackrel{\cong}{\to} A, \quad \text{ for all } A\in{\rm Ob}\,\mathcal{C},
\]
a unit for $\square$, is called a {\em tensor category} or {\em monoidal category}, if the {\em pentagon axiom}
\[
\xymatrix{
 & (A\square B)\square(C\square D)\ar[dr]^{\alpha_{A\square B,C,C}} & \\
A\square \big(B \square (C \square D)\big) \quad\ar[ur]^{\alpha_{A,B,C\square D}}\ar[d]_{{\rm id}_A\square\alpha_{A,B,C}} & & \quad\big((A\square B) \square C\big) \square D \\
  \qquad A\square \big((B \square C) \square D\big)\ar[rr]_{\alpha_{A,B\square C.D}} & & \big(A\square (B \square C)\big) \square D \ar[u]_{\alpha_{A,B,C}\square{\rm id}_D} \qquad
}
\]
and the {\rm triangle axiom}
\[
\xymatrix{
A\square (E\square C)\ar[rr]^{\alpha_{A,E,C}}\ar[dr]_{{\rm id}_A\square \lambda_C} & & (A\square E)C\ar[dl]^{\rho_A\square{\rm id}_C} \\
& A\square C & \\
}
\]
are satisfied for all objects $A,B,C,D$ of $\mathcal{C}$.
\end{definition}

If a category has products or coproducts for all finite sets of objects, then the universal property guarantees the existence of the isomorphisms $\alpha$, $\lambda$, and $\rho$ that turn it into a tensor category.

In order to define a notion of independence we need less than a
(co-) product, but a little bit more than a
tensor product. What we need are inclusions or projections that allow
us to view the objects $A$, $B$ as subsystems of their product $A\square B$.

\begin{definition}\label{def-tensor-cat-with}
A {\em tensor category with projections} $(\mathcal{C},\square,\pi)$
is a tensor category $(\mathcal{C},\square)$ equipped with two natural transformations $\pi_1:\square \to P_1$ and $\pi_2:\square \to P_2$, where the bifunctors $P_1,P_2:\mathcal{C}\times\mathcal{C}\to\mathcal{C}$ are defined by $P_1(B_1,B_2)=B_1$, $P_2(B_1,B_2)=B_2$, on pairs of objects $B_1,B_2$ of $\mathcal{C}$, and similarly on pairs of morphisms. In other words, for any pair of objects $B_1,B_2$ there exist two morphisms $\pi_{B_1}:B_1\square B_2\to B_1$, $\pi_{B_2}:B_1\square B_2\to B_2$, such that for any pair of morphisms $f_1:A_1\to B_1$, $f_2:A_2\to B_2$, the following diagram commutes,
\[
\xymatrix{
A_1 \ar[d]_{f_1}& A_1\square A_2 \ar[d]|-{f_1\square f_2} \ar[l]_{\pi_{A_1}} \ar[r]^{\pi_{A_2}} & A_2\ar[d]^{f_2} \\
B_1 & B_1\square B_2\ar[l]^{\pi_{B_1}} \ar[r]_{\pi_{B_2}} & B_2 .
}
\]

Similarly, a {\em tensor product with inclusions}
$(\mathcal{C},\square,i)$ is a tensor category $(\mathcal{C},\square)$
equipped with two natural transformations $i_1:P_1\to\square$ and
$i_2:P_2\to\square$,  i.e.\ for any pair of objects $B_1,B_2$ there exist two morphisms $i_{B_1}:B_1\to B_1\square B_2$, $i_{B_2}:B_2\to B_1\square B_2$, such that for any pair of morphisms $f_1:A_1\to B_1$, $f_2:A_2\to B_2$, the following diagram commutes,
\[
\xymatrix{
A_1 \ar[r]_{i_{A_1}} \ar[d]_{f_1}& A_1\square A_2 \ar[d]|-{f_1\square f_2} & A_2\ar[d]^{f_2} \ar[l]^{i_{A_2}} \\
B_1 \ar[r]^{i_{B_1}}& B_1\square B_2 & B_2\ar[l]_{i_{B_2}} .
}
\]
\end{definition}

In  a tensor category with projections or with inclusions we can define a notion of independence for morphisms.

\begin{definition}\label{def-independence}
Let $(\mathcal{C},\square,\pi)$ be a tensor category with
projections. Two morphism $f_1:A\to B_1$ and $f_2:A\to B_2$ with the same
source $A$ are called
{\em independent}\index{independent} (with respect to $\square$), if there exists a morphism $h:A\to B_1\square B_2$ such that the
diagram
\begin{equation}\label{diag-indep1}
\xymatrix{
& A \ar[dl]_{f_1} \ar@{.>}[d]|-h\ar[dr]^{f_2} & \\
B_1 & B_1\square B_2 \ar[l]^{\pi_{B_1}} \ar[r]_{\pi_{B_2}} & B_2
}
\end{equation}
commutes.

In a tensor category with inclusions $(\mathcal{C},\square,i)$, two morphisms $f_1:A_1\to B$ and $f_2:A_2\to B$ with the same
target $B$ are called independent, if there exists a morphism $h:A_1\square
A_2\to B$ such that the diagram
\begin{equation}\label{diag-indep2}
\xymatrix{
& A  & \\
B_1 \ar[ur]^{f_1}\ar[r]_{i_{B_1}}& B_1\square B_2 \ar@{.>}[u]|-h & B_2\ar[l]^{i_{B_2}}\ar[ul]_{f_2}
}
\end{equation}
commutes.
\end{definition}
This definition can be extended in the obvious way to arbitrary sets of morphisms.

If $\square$ is actually a product (or coproduct, resp.), then the universal property in Definition \ref{def-product} implies that for all pairs of morphisms with the same source (or target,
resp.) there exists even a unique morphism
that makes diagram \eqref{diag-indep1} (or \eqref{diag-indep2}, resp.)
commuting. Therefore in that case
all pairs of morphism with the same source (or target, resp.) are
independent.

We will now consider several examples. We will show that for the category of classical probability spaces we recover usual stochastic independence, if we take the product of probability spaces, cf.\ Proposition \ref{prop-stoch-independence}.

\subsection{Example: Independence in the Category of Classical Probability Spaces}\label{example-class-indep}

The category $\eufrak{Meas}$ of measurable spaces consists of pairs $(\Omega,\mathcal{F})$, where $\Omega$ is a set and $\mathcal{F}\subseteq\mathcal{P}(\Omega)$ a $\sigma$-algebra. The morphisms are the measurable maps. This category has a product,
\[
(\Omega_1,\mathcal{F}_1)\,\Pi\,(\Omega_2,\mathcal{F}_2)=(\Omega_1\times\Omega_2,\mathcal{F}_1\otimes\mathcal{F}_2)
\]
where $\Omega_1\times\Omega_2$ is the Cartesian product of $\Omega_1$ and $\Omega_2$, and $\mathcal{F}_1\otimes\mathcal{F}_2$ is the smallest $\sigma$-algebra on $\Omega_1\times\Omega_2$ such that the canonical projections $p_1:\Omega_1\times\Omega_2\to\Omega_1$ and $p_2:\Omega_1\times\Omega_2\to\Omega_2$ are measurable.

The category of probability spaces $\eufrak{Prob}$ has as objects triples $(\Omega,\mathcal{F},P)$ where $(\Omega,\mathcal{F})$ is a measurable space and $P$ a probability measure on $(\Omega,\mathcal{F})$. A morphism $X:(\Omega_1,\mathcal{F}_1,P_1)\to (\Omega_1,\mathcal{F}_2,P_2)$ is a measurable map $X:(\Omega_1,\mathcal{F}_1)\to (\Omega_1,\mathcal{F}_2)$ such that
\[
P_1\circ X^{-1} = P_2.
\]
This means that a random variable $X:(\Omega,\mathcal{F})\to(E,\mathcal{E})$ automatically becomes a morphism, if we equip $(E,\mathcal{E})$ with the measure
\[
P_X=P\circ X^{-1}
\]
induced by $X$.

This category does not have universal products. But one can check that the product of measures turns $\eufrak{Prob}$ into a tensor category,
\[
(\Omega_1,\mathcal{F}_1,P_1)\otimes(\Omega_2,\mathcal{F}_2,P_2)=(\Omega_1\times\Omega_2,\mathcal{F}_1\otimes\mathcal{F}_2,P_1\otimes P_2),
\]
where $P_1\otimes P_2$ is determined by
\[
(P_1\otimes P_2)(M_1\times M_2) =P_1(M_1)P_2(M_2),
\]
for all $M_1\in\mathcal{F}_1$, $M_2\in\mathcal{F}_2$. It is even a tensor category with projections in the sense of Definition \ref{def-tensor-cat-with} with the canonical projections $p_1:(\Omega_1\times\Omega_2,\mathcal{F}_1\otimes\mathcal{F}_2,P_1\otimes P_2)\to (\Omega_1,\mathcal{F}_1,P_1)$, $p_2:(\Omega_1\times\Omega_2,\mathcal{F}_1\otimes\mathcal{F}_2,P_1\otimes P_2)\to (\Omega_2,\mathcal{F}_2,P_2)$ given by $p_1\big((\omega_1,\omega_2)\big)=\omega_1$, $p_2\big((\omega_1,\omega_2)\big)=\omega_2$ for $\omega_1\in\Omega_1$, $\omega_2\in\Omega_2$.

The notion of independence associated to this tensor product with projections is exactly the one used in probability.

\begin{proposition}\label{prop-stoch-independence}
Two random variables $X_1:(\Omega,\mathcal{F},P)\to(E_1,\mathcal{E}_1)$ and $X_2:(\Omega,\mathcal{F},P)\to(E_2,\mathcal{E}_2)$, defined on the same probability space $(\Omega,\mathcal{F},P)$ and with values in measurable spaces $(E_1,\mathcal{E}_1)$ and $(E_2,\mathcal{E}_2)$, are stochastically independent, if and only if they are independent in the sense of Definition \ref{def-independence} as morphisms $X_1:(\Omega,\mathcal{F},P)\to(E_1,\mathcal{E}_1,P_{X_1})$ and $X_2:(\Omega,\mathcal{F},P)\to(E_2,\mathcal{E}_2,P_{X_2})$ of the tensor category with projections $(\eufrak{Prob},\otimes,p)$.
\end{proposition}
\begin{proof}
Assume that $X_1$ and $X_2$ are stochastically independent. We have to find a morphism $h:(\Omega,\mathcal{F},P)\to(E_1\times E_2,\mathcal{E}_1\otimes\mathcal{E}_2,P_{X_1}\otimes P_{X_2})$ such that the diagram
\[
\xymatrix{
& (\Omega,\mathcal{F},P) \ar[dl]_{X_1} \ar@{.>}[d]|-h\ar[dr]^{X_2} & \\
(E_1,\mathcal{E}_1,P_{X_1}) &  (E_1\times E_2,\mathcal{E}_1\otimes\mathcal{E}_2,P_{X_1}\otimes P_{X_2})\ar[l]^(.65){p_{E_1}} \ar[r]_(.65){p_{E_2}} & (E_2,\mathcal{E}_2,P_{X_2})
}
\]
commutes. The only possible candidate is $h(\omega)=\big(X_1(\omega),X_2(\omega)\big)$ for all $\omega\in\Omega$, the unique map that completes this diagram in the category of measurable spaces and that exists due to the universal property of the product of measurable spaces. This is a morphism in $\eufrak{Prob}$, because we have
\begin{eqnarray*}
P\big(h^{-1}(M_1\times M_2)\big)&=&P\big(X_1^{-1}(M_1)\cap X_2^{-1}(M_2)\big) = P\big(X_1^{-1}(M_1)\big)P\big(X_2^{-1}(M_2)\big)\\
&=& P_{X_1}(M_1)P_{X_2}(M_2)=(P_{X_1}\otimes P_{X_2})(M_1\times M_2)
\end{eqnarray*}
for all $M_1\in\mathcal{E}_1$, $M_2\in\mathcal{E}_2$, and therefore
\[
P\circ h^{-1}=P_{X_1}\otimes P_{X_2}.
\]
Conversely, if $X_1$ and $X_2$ are independent in the sense of Definition \ref{def-independence}, then the morphism that makes the diagram commuting has to be again $h:\omega\mapsto \big(X_1(\omega),X_2(\omega)\big)$. This implies
\[
P_{(X_1,X_2)}=P\circ h^{-1} = P_{X_1}\otimes P_{X_2}
\]
and therefore
\[
P\big(X_1^{-1}(M_1)\cap X_2^{-1}(M_2)\big) =  P\big(X_1^{-1}(M_1)\big)P\big(X_2^{-1}(M_2)\big)
\]
for all $M_1\in\mathcal{E}_1$, $M_2\in\mathcal{E}_2$.
\end{proof}

\subsection{Example: Tensor Independence in the Category of Algebraic Probability Spaces}\label{tensor indep}

By the category of algebraic probability $\eufrak{AlgProb}$ spaces we denote the category of associative unital algebras over $\mathbb{C}$ equipped with a unital linear functional. A morphism $j:(\mathcal{A}_1,\varphi_1)\to(\mathcal{A}_2,\varphi_2)$ is a quantum random variable, i.e.\ an algebra homomorphism $j:\mathcal{A}_1\to\mathcal{A}_2$ that preserves the unit and the functional, i.e.\ $j(\mathbf{1}_{\mathcal{A}_1})=\mathbf{1}_{\mathcal{A}_2}$ and $\varphi_2\circ j=\varphi_1$.

The tensor product we will consider on this category is just the usual tensor product $(\mathcal{A}_1\otimes\mathcal{A}_2,\varphi_1\otimes\varphi_2)$, i.e.\ the algebra structure of $\mathcal{A}_1\otimes\mathcal{A}_2$ is defined by
\begin{eqnarray*}
\mathbf{1}_{\mathcal{A}_1\otimes\mathcal{A}_2} &=& \mathbf{1}_{\mathcal{A}_1}\otimes\mathbf{1}_{\mathcal{A}_2}, \\
(a_1\otimes a_2)(b_1\otimes b_1) &=& a_1b_2 \otimes a_2b_2,
\end{eqnarray*}
and the new functional is defined by
\[
(\varphi_1\otimes\varphi_2)(a_1\otimes a_2) = \varphi_1(a_1)\varphi_2(a_2),
\]
for all $a_1,b_1\in\mathcal{A}_1$, $a_2,b_2\in\mathcal{A}_2$.

This becomes a tensor category with inclusions with the inclusions defined by
\begin{eqnarray*}
i_{\mathcal{A}_1}(a_1) &=& a_1\otimes \mathbf{1}_{\mathcal{A}_2}, \\
i_{\mathcal{A}_2}(a_2) &=& \mathbf{1}_{\mathcal{A}_1}\otimes a_2,
\end{eqnarray*}
for $a_1\in\mathcal{A}_1$, $a_2\in\mathcal{A}_2$.

One gets the category of $*$-algebraic probability spaces, if one assumes that the underlying algebras have an involution and the functional are states, i.e.\ are also positive. Then an involution is defined on $\mathcal{A}_1\otimes\mathcal{A}_2$ by $(a_1\otimes a_2)^*= a_1^*\otimes a_2^*$ and $\varphi_1\otimes\varphi_2$ is again a state.

The notion of independence associated to this tensor product with inclusions by Definition \ref{def-independence} is the usual notion of {\em Bose} or {\em tensor independence} used in quantum probability, e.g., by Hudson and Parthasarathy.

\begin{proposition}
Two quantum random variables $j_1:(\mathcal{B}_1,\psi_1)\to(\mathcal{A},\varphi)$ and $j_2:(\mathcal{B}_2,\psi_2)\to(\mathcal{A},\varphi)$, defined on algebraic probability spaces $(\mathcal{B}_1,\psi_1),(\mathcal{B}_2,\psi_2)$ and with values in the same algebraic probability space $(\mathcal{A},\varphi)$ are independent if and only if the following two conditions are satisfied.
\begin{description}
\item[(i)]
The images of $j_1$ and $j_2$ commute, i.e.
\[
\big[j_1(a_1),j_2(a_2)\big] = 0,
\]
for all $a_1\in\mathcal{A}_1$, $a_2\in\mathcal{A}_2$.
\item[(ii)]
$\varphi$ satisfies the factorization property
\[
\varphi\big(j_1(a_1)j_2(a_2)\big) = \varphi\big(j_1(a_1)\big)\varphi\big(j_2(a_2)\big),
\]
for all $a_1\in\mathcal{A}_1$, $a_2\in\mathcal{A}_2$.
\end{description}
\end{proposition}
We will not prove this Proposition since it can be obtained as a special case of Proposition \ref{prop-fermi-indep}, if we equip the algebras with the trivial $\mathbb{Z}_2$-grading $\mathcal{A}^{(0)}=\mathcal{A}$, $\mathcal{A}^{(1)}=\{0\}$.

\subsection{Example: Fermi Independence}\label{fermi-indep}

Let us now consider the category of $\mathbb{Z}_2$-graded algebraic probability spaces $\mathbb{Z}_2$-$\eufrak{AlgProb}$. The objects are pairs $(\mathcal{A},\varphi)$ consisting of a $\mathbb{Z}_2$-graded unital algebra $\mathcal{A}=\mathcal{A}^{(0)}\oplus\mathcal{A}^{(1)}$ and an even unital functional $\varphi$, i.e.\ $\varphi|_{\mathcal{A}^{(1)}}=0$. The morphisms are random variables that don't change the degree, i.e., for $j:(\mathcal{A}_1,\varphi_1)\to(\mathcal{A}_2,\varphi_2)$, we have
\[
j(\mathcal{A}_1^{(0)})\subseteq\mathcal{A}_2^{(0)} \quad \mbox{ and }\quad j(\mathcal{A}_1^{(1)})\subseteq\mathcal{A}_2^{(1)}.
\]
The tensor product $(\mathcal{A}_1\otimes_{\mathbb{Z}_2}\mathcal{A}_2,\varphi_1\otimes\varphi_2)= (\mathcal{A}_1,\varphi_1)\otimes_{\mathbb{Z}_2}(\mathcal{A}_2,\varphi_2)$ is defined as follows. The algebra $\mathcal{A}_1\otimes_{\mathbb{Z}_2}\mathcal{A}_2$ is the graded tensor product of $\mathcal{A}_1$ and $\mathcal{A}_2$, i.e.\ $(\mathcal{A}_1\otimes_{\mathbb{Z}_2}\mathcal{A}_2)^{(0)}=\mathcal{A}^{(0)}_1\otimes\mathcal{A}^{(0)}_2\oplus \mathcal{A}^{(1)}_1\otimes\mathcal{A}^{(1)}_2$,  $(\mathcal{A}_1\otimes_{\mathbb{Z}_2}\mathcal{A}_2)^{(1)}=\mathcal{A}^{(1)}_1\otimes\mathcal{A}^{(0)}_2\oplus \mathcal{A}^{(0)}_1\otimes\mathcal{A}^{(1)}_2$, with the algebra structure given by
\begin{eqnarray*}
\mathbf{1}_{\mathcal{A}_1\otimes_{\mathbb{Z}_2}\mathcal{A}_2} &=& \mathbf{1}_{\mathcal{A}_1} \otimes\mathbf{1}_{\mathcal{A}_2}, \\
(a_1\otimes b_1) \cdot (a_2\otimes b_2) &=& (-1)^{\deg b_1\deg a_2} a_1b_1\otimes a_2b_2, 
\end{eqnarray*}
for all homogeneous elements $a_1,b_1\in\mathcal{A}_1$, $a_2,b_2\in\mathcal{A}_2$. The functional $\varphi_1\otimes\varphi_2$ is simply the tensor product, i.e.\ $(\varphi_1\otimes\varphi_2)(a_1\otimes a_2)=\varphi_1(a_1)\otimes \varphi_2(a_2)$ for all $a_1\in\mathcal{A}_1$, $a_2\in\mathcal{A}_2$. It is easy to see that $\varphi_1\otimes\varphi_2$ is again even, if $\varphi_1$ and $\varphi_2$ are even. The inclusions $i_1:(\mathcal{A}_1,\varphi_1)\to(\mathcal{A}_1\otimes_{\mathbb{Z}_2}\mathcal{A}_2,\varphi_1\otimes\varphi_2)$ and  $i_2:(\mathcal{A}_2,\varphi_2)\to(\mathcal{A}_1\otimes_{\mathbb{Z}_2}\mathcal{A}_2,\varphi_1\otimes\varphi_2)$ are defined by
\[
i_1(a_1) = a_1\otimes \mathbf{1}_{\mathcal{A}_2}\quad\mbox{ and } i_2(a_2)= \mathbf{1}_{\mathcal{A}_1}\otimes a_2,
\]
for $a_1\in\mathcal{A}_1$, $a_2\in\mathcal{A}_2$.

If the underlying algebras are assumed to have an involution and the functionals to be states, then the involution on the $\mathbb{Z}_2$-graded tensor product is defined by $(a_1\otimes a_2)^* = (-1)^{\deg a_1 \deg a_2} a_1^*\otimes a_2^*$, this gives the category of $\mathbb{Z}_2$-graded $*$-algebraic probability spaces.

The notion of independence associated to this tensor category with inclusions is called {\em Fermi independence} or {\em anti-symmetric independence}.

\begin{proposition}\label{prop-fermi-indep}
Two random variables
$j_1:(\mathcal{B}_1,\psi_1)\to(\mathcal{A},\varphi)$ and
$j_2:(\mathcal{B}_2,\psi_2)\to(\mathcal{A},\varphi)$, defined on two
$\mathbb{Z}_2$-graded algebraic probability spaces
$(\mathcal{B}_1,\psi_1)$, $(\mathcal{B}_2,\psi_2)$ and with values in the same $\mathbb{Z}_2$-algebraic probability space $(\mathcal{A},\varphi)$ are independent if and only if the following two conditions are satisfied.
\begin{description}
\item[(i)]
The images of $j_1$ and $j_2$ satisfy the commutation relations
\[
j_2(a_2)j_1(a_1) = (-1)^{\deg a_1 \deg a_2}j_1(a_1)j_2(a_2)
\]
for all homogeneous elements $a_1\in\mathcal{A}_1$, $a_2\in\mathcal{A}_2$.
\item[(ii)]
$\varphi$ satisfies the factorization property
\[
\varphi\big(j_1(a_1)j_2(a_2)\big) = \varphi\big(j_1(a_1)\big)\varphi\big(j_2(a_2)\big),
\]
for all $a_1\in\mathcal{A}_1$, $a_2\in\mathcal{A}_2$.
\end{description}
\end{proposition}
\begin{proof}
The proof is similar to that of Proposition \ref{prop-stoch-independence}, we will only outline it. It is clear that the morphism $h:(\mathcal{B}_1,\psi_1)\otimes_{\mathbb{Z}_2}(\mathcal{B}_2,\psi_2)\to(\mathcal{A},\varphi)$ that makes the diagram in Definition \ref{def-independence} commuting, has to act on elements of $\mathcal{B}_1\otimes \mathbf{1}_{\mathcal{B}_2}$ and $\mathbf{1}_{\mathcal{B}_1}\otimes\mathcal{B}_2$ as
\[
h(b_1\otimes \mathbf{1}_{\mathcal{B}_2}) = j_1(b_1) \quad\mbox{ and } \quad h(\mathbf{1}_{\mathcal{B}_1}\otimes b_2)=j_2(b_2).
\]
This extends to a homomorphism from $(\mathcal{B}_1,\psi_1)\otimes_{\mathbb{Z}_2}(\mathcal{B}_2,\psi_2)$ to $(\mathcal{A},\varphi)$, if and only if the commutation relations are satisfied. And the resulting homomorphism is a quantum random variable, i.e.\ satisfies $\varphi\circ h=\psi_1\otimes\psi_2$, if and only if the factorization property is satisfied.
\end{proof}

\section{Reduction of Independences}\label{reduction}

In this Section we will study the relations between different notions of independence. Let us first recall the definition of a tensor functor.

\begin{definition}(see, e.g., Section XI.2 in MacLane\cite{maclane98})
Let $(\mathcal{C},\square)$ and $(\mathcal{C}',\square')$ be two tensor categories.
A {\em cotensor functor} or {\em comonoidal functor} $F:(\mathcal{C},\square)\to(\mathcal{C}',\square')$ is an ordinary functor $F:\mathcal{C}\to\mathcal{C}'$ equipped with a morphism $F_0: F(E_\mathcal{C}) \to E_{\mathcal{C}'}$ and a natural transformation $F_2:F(\,\cdot\, \square \,\cdot\,) \to  F(\,\cdot\,) \square' F(\,\cdot\,)$, i.e.\ morphisms $F_2(A,B):F(A \square B) \to  F(A) \square' F(B)$ for all $A,B\in{\rm Ob}\,\mathcal{C}$ that are natural in $A$ and $B$, such that the diagrams
\begin{equation}
\xymatrix{
F\big(A\square (B\square C)\big)\ar[rr]^{F(\alpha_{A,B,C})} \ar[d]_{F_2(A,B\square C)} && F\big((A\square B)\square C\big)\ar[d]^{F_2(A\square B,C)} \\
F(A)\square' F(B\square C) \ar[d]_{{\rm id}_{F(A)}\square'F_2(B,C)} && F(A\square B)\square' F(C)\ar[d]^{F_2(A,B)\square'{\rm id}_{F(C)}} \\
F(A)\square'\big(F(B)\square'F(C)\big) \ar[rr]_{\alpha'_{F(A),F(B),F(C)}}&& \big(F(A)\square' F(B)\big)\square' F(C)
}
\end{equation}
\begin{equation}
\xymatrix{
F(B\square E_\mathcal{C}) \ar[rr]^{F_2(B,E_\mathcal{C})} \ar[d]_{F(\rho_B)} && F(B)\square' F(E_\mathcal{C})\ar[d]^{{\rm id}_B\square' F_0} \\
F(B) && F(B)\square'E_{\mathcal{C}'} \ar[ll]^{\rho_{F(B)}}
}
\end{equation}
\begin{equation}
\xymatrix{
F(E_\mathcal{C}\square B) \ar[rr]^{F_2(E_\mathcal{C},B)} \ar[d]_{F(\lambda_B)} && F(E_\mathcal{C})\square'F(B) \ar[d]^{F_0\square'{\rm id}_B} \\
F(B) && E_{\mathcal{C}'} \square'F(B) \ar[ll]^{\lambda_{F(B)}}
}
\end{equation}
commute for all $A,B,C\in{\rm Ob}\,\mathcal{C}$.

\end{definition}
We have reversed the direction of $F_0$ and $F_2$ in our definition. In the case of a strong tensor functor, i.e.\ when all the morphisms are isomorphisms, our definition of a cotensor functor is equivalent to the usual definition of a tensor functor as, e.g., in MacLane\cite{maclane98}.

The conditions are exactly what we need to get morphisms $F_n(A_1,\ldots,A_n):F(A_1\square\cdots\square A_n)\to F(A_1)\square'\cdots\square F(A_n)$ for all finite sets $\{A_1,\ldots,A_n\}$ of objects of $\mathcal{C}$ such that, up to these morphisms, the functor $F:(\mathcal{C},\square)\to(\mathcal{C}',\square')$ is a homomorphism.

For a reduction of independences we need a little bit more than a cotensor functor.

\begin{definition}\label{def-reduction}
Let $(\mathcal{C},\square,i)$ and $(\mathcal{C}',\square',i')$ be two tensor categories with inclusions and assume that $\mathcal{C}$ is a subcategory of $\mathcal{C}'$. A {\em reduction} $(F,J)$ of the tensor product $\square$ to the tensor product $\square'$ is a cotensor functor $F:(\mathcal{C},\square)\to(\mathcal{C}',\square')$ and a natural transformation $J:{\rm id}_\mathcal{C}\to F$, i.e.\ morphisms $j_A:A\to F(A)$ in $\mathcal{C}'$ for all objects $A\in{\rm Ob}\,\mathcal{C}$ such that the diagram
\[
\xymatrix{
A \ar[d]_f \ar[r]^{J_A} & F(A)\ar[d]^{F(f)} \\
B \ar[r]_{J_B} & F(B)
}
\]
commutes for all morphisms $f:A\to B$ in $\mathcal{C}$.
\end{definition}
Such a reduction provides us with a system of inclusions
$J_n(A_1,\ldots,A_n)=F_n(A_1,\ldots,A_n)\circ J_{A_1\square\cdots
  \square A_n}:A_1\square\cdots \square A_n\to
F(A_1)\square'\cdots\square'F(A_n)$ with $J_1(A)=J_A$ that satisfies,
e.g., $J_{n+m}(A_1,\ldots,A_{n+m}) = F_2\big(F(A_1)\square' \cdots
\square'F(A_n),$ $F(A_{n+1})\square' \cdots \square'F(A_{n+m})\big) \circ \big(J_n(A_1, \ldots ,A_n)\square J_m(A_{n+1}, \ldots ,A_{n+m})\big)$ for all $n,m\in\mathbb{N}$ and $A_1,\ldots,A_{n+m}\in{\rm Ob}\,\mathcal{C}$.

A reduction between two tensor categories with projections would consist of a tensor functor $F$ and a natural transformation $P:F\to{\rm id}$.

We have to extend our definition slightly. In our applications $\mathcal{C}$ will often not be a subcategory of $\mathcal{C}'$, but we have, e.g., a forgetful functor $U$ from $\mathcal{C}$ to $\mathcal{C}'$ that ``forgets'' an additional structure that $\mathcal{C}$ has. An example for this situation is the reduction of Fermi independence to tensor independence in following subsection. Here we have to forget the $\mathbb{Z}_2$-grading of the objects of $\mathbb{Z}_2$-$\eufrak{AlgProb}$ to get objects of $\eufrak{AlgProb}$. In this situation a reduction of the tensor product with inclusions $\square$ to the tensor product with inclusions $\square'$ is a tensor function $F$ from $(\mathcal{C},\square)$ to $(\mathcal{C}',\square')$ and a natural transformation $J:U\to F$.

\subsection{Example: Bosonization of Fermi Independence}\label{reduction-fermi-bose}

We will now define the bosonization of Fermi independence as a reduction from $(\eufrak{AlgProb},\otimes,i)$ to $(\mathbb{Z}_2$-$\eufrak{AlgProb},\otimes_{\mathbb{Z}_2},i)$. We will need the group algebra $\mathbb{C}\mathbb{Z}_2$ of $\mathbb{Z}_2$ and the linear functional $\varepsilon:\mathbb{C}\mathbb{Z}_2\to\mathbb{C}$ that arises as the linear extension of the trivial representation of $\mathbb{Z}_2$, i.e.
\[
\varepsilon(\mathbf{1})=\varepsilon(g)=1,
\]
if we denote the even element of $\mathbb{Z}_2$ by $\mathbf{1}$ and the odd element by $g$.

The underlying functor $F:\mathbb{Z}_2$-$\eufrak{AlgProb}\to\eufrak{AlgProb}$ is given by
\[
F:
\begin{array}{lcl}
{\rm Ob}\,\mathbb{Z}_2\mbox{-}\eufrak{AlgProb}\ni(\mathcal{A},\varphi) & \mapsto & (\mathcal{A}\otimes_{\mathbb{Z}_2}\mathbb{C}\mathbb{Z}_2, \varphi\otimes\varepsilon)\in{\rm Ob}\,\eufrak{AlgProb}, \\[1mm]
{\rm Mor}\,\mathbb{Z}_2\mbox{-}\eufrak{AlgProb}\ni f & \mapsto &  f\otimes {\rm id}_{\mathbb{C}\mathbb{Z}_2} \in{\rm Mor}\,\eufrak{AlgProb}.
\end{array}
\]

The unit element in both tensor categories is the one-dimensional unital algebra $\mathbb{C}\mathbf{1}$ with the unique unital functional on it. Therefore $F_0$ has to be a morphism from $F(\mathbb{C}\mathbf{1})\cong \mathbf{C}\mathbb{Z}_2$ to $\mathbb{C}\mathbf{1}$. It is defined by $F_0(\mathbf{1})=F_0(g)=\mathbf{1}$.

The morphism $F_2(\mathcal{A}_1,\mathcal{A}_2)$ has to go from $F(\mathcal{A}\otimes_{\mathbb{Z}_2}\mathcal{B})=(\mathcal{A}\otimes_{\mathbb{Z}_2}\mathcal{B})\otimes \mathbb{C}\mathbb{Z}_2$ to $F(\mathcal{A})\otimes F(\mathcal{B})=(\mathcal{A}\otimes_{\mathbb{Z}_2}\mathbf{C}\mathbb{Z}_2)\otimes(\mathcal{B}\otimes_{\mathbb{Z}_2}\mathbf{C}\mathbb{Z}_2)$. It is defined by
\[
a\otimes b \otimes \mathbf{1} \mapsto
\left\{
\begin{array}{lcl}
(a\otimes \mathbf{1}) \otimes (b\otimes \mathbf{1}) & \mbox{ if } & b\mbox{ is even}, \\
(a\otimes g) \otimes (b\otimes \mathbf{1} & \mbox{ if } & b\mbox{ is odd},
\end{array}
\right.
\]
and
\[
a\otimes b \otimes g \mapsto
\left\{
\begin{array}{lcl}
(a\otimes g) \otimes (b\otimes g) & \mbox{ if } & b\mbox{ is even}, \\
(a\otimes \mathbf{1}) \otimes (b\otimes g) & \mbox{ if } & b\mbox{ is odd},
\end{array}
\right.
\]
for $a\in\mathcal{A}$ and homogeneous $b\in\mathcal{B}$.

Finally, the inclusion $J_\mathcal{A}:\mathcal{A}\to \mathcal{A}\otimes_{\mathbb{Z}_2}\mathbf{C}\mathbb{Z}_2$ is defined by
\[
J_\mathcal{A}(a) = a\otimes\mathbf{1}
\]
for all $a\in\mathcal{A}$.

In this way we get inclusions $J_n=J_n(\mathcal{A}_1,\ldots,\mathcal{A}_n)=F_n(\mathcal{A}_1,\ldots,\mathcal{A}_n)\circ J_{\mathcal{A}_1\otimes_{\mathbb{Z}_2}\ldots\otimes_{\mathbb{Z}_2}\mathcal{A}_n}$ of the graded tensor product $\mathcal{A}_1\otimes_{\mathbb{Z}_2}\cdots\otimes_{\mathbb{Z}_2}\mathcal{A}_n$ into the usual tensor product $(\mathcal{A}_1\otimes_{\mathbb{Z}_2}\mathbb{C}\mathbb{Z}_2)\otimes\cdots\otimes (\mathcal{A}_n\otimes_{\mathbb{Z}_2}\mathbb{C}\mathbb{Z}_2)$ which respect the states and allow to reduce all calculations involving the graded tensor product to calculations involving the usual tensor product on the bigger algebras $F(\mathcal{A}_1)=\mathcal{A}_1\otimes_{\mathbb{Z}_2}\mathbb{C}\mathbb{Z}_2,\ldots,F(\mathcal{A}_n)=\mathcal{A}_n\otimes_{\mathbb{Z}_2}\mathbb{C}\mathbb{Z}_2$. These inclusions are determined by
\begin{eqnarray*}
J_n(\underbrace{\mathbf{1}\otimes\cdots\otimes\mathbf{1}}\otimes a \otimes\underbrace{\mathbf{1}\otimes\cdots\otimes\mathbf{1}}) &=& \underbrace{\tilde{g}\otimes\cdots\otimes\tilde{g}}\otimes \tilde{a}\otimes\underbrace{\tilde{\mathbf{1}}\otimes\cdots\otimes\tilde{\mathbf{1}}}, \\
\hspace*{11mm} k-1 \mbox{ times} \hspace{8mm} n-k \mbox{ times} && k-1 \mbox{ times} \hspace{8mm} n-k \mbox{ times}
\end{eqnarray*}
for $a\in \mathcal{A}_k$, $1\le k\le n$, where we used the abbreviations
\[
\tilde{g}=\mathbf{1}\otimes g, \qquad \tilde{a}=a\otimes\mathbf{1}, \qquad \tilde{\mathbf{1}}=\mathbf{1}\otimes\mathbf{1}.
\]

\section{Forgetful Functors, Coproducts, and Semi-universal Products}

We are mainly interested in different categories of algebraic
probability spaces. There objects are pairs consisting of an algebra
$\mathcal{A}$ and a linear functional $\varphi$ on
$\mathcal{A}$. Typically, the algebra has some additional structure,
e.g., an involution, a unit, a grading, or a topology (it can be, e.g., a von Neumann algebra or a $C^*$-algebra), and the functional behaves nicely with respect to this additional structure, i.e., it is positive, unital, respects the grading, continuous, or normal. The morphisms are algebra homomorphisms, which leave the linear functional invariant, i.e., $j:(\mathcal{A},\varphi)\to (\mathcal{B},\psi)$ satisfies
\[
\varphi=\psi\circ j
\]
and behave also nicely w.r.t.\ to additional structure, i.e., they can be
required to be $*$-algebra homomorphisms, map the unit of
$\mathcal{A}$ to the unit of $\mathcal{B}$, respect the grading, etc. We have already seen one example in Subsection \ref{fermi-indep}.

The tensor product then has to specify a new algebra with a linear functional
and inclusions for every pair of of algebraic probability spaces. If the category of algebras obtained from our algebraic
probability space by forgetting the linear functional has a coproduct, then it
is sufficient to consider the case where the new algebra is the coproduct of
the two algebras. 

\begin{proposition}\label{prop-forget}
Let $(\mathcal{C},\square,i)$ be a tensor category with inclusions and $F:\mathcal{C}\to\mathcal{D}$ a functor from $\mathcal{C}$ into another category $\mathcal{D}$ which has a coproduct $\coprod$ and an initial object $E_\mathcal{D}$. Then $F$ is a tensor functor. The morphisms $F_2(A,B):F(A)\coprod F(B)\to F(A\square B)$ and $F_0:E_\mathcal{D}\to F(E)$ are those guaranteed by the universal property of the coproduct and the initial object, i.e.\ $F_0:E_\mathcal{D}\to F(E)$ is the unique morphism from $E_\mathcal{D}$ to $F(E)$ and $F_2(A,B)$ is the unique morphism that makes the diagram
\[
\xymatrix{
F(A) \ar[r]^(0.4){F(i_A)}\ar[dr]_{i_{F(A)}} & F(A\square B) & \ar[l]_(0.4){F(i_B)} \ar[dl]^{i_{F(B)}} F(B) \\
 & F(A)\coprod F(B) \ar@{.>}[u]|-{F_2(A,B)} &
}
\]
commuting.
\end{proposition}
\begin{proof}
Using the universal property of the coproduct and the definition of $F_2$, one shows that the triangles containing the $F(A)$ in the center of the diagram
\[
\xymatrix{
F(A)\coprod \big(F(B)\coprod F(C)\big) \ar[rr]^{\alpha_{F(A),F(B),F(C)}} \ar[d]_{{\rm id}_{F(A)}\coprod F_2(B,C)} && \big(F(A)\coprod F(B)\big)\coprod F(C)\ar[d]^{F_2(A,B)\coprod{\rm id}_{F(C)}} \\
F(A)\coprod F(B\square C)\ar[d]_{F_2(A,B\square C)} & F(A) \ar[ul]|-{i_{F(A)}} \ar[l]|-{i_{F(A)}} \ar[dl]|-{F(i_A)} \ar[ur]|-{i_{F(A)}} \ar[r] \ar[rd]|-{F(i_A)} & F(A\square B)\coprod F(C)\ar[d]^{F_2(A\square B,C)} \\
F\big(A\square(B\square C)\big)\ar[rr]_{F(\alpha_{A,B,C})} & & F\big((A\square B)\square C\big)
}
\]
commute (where the morphism from $F(A)$ to $ F(A\square B)\coprod F(C)$ is $F(i_A)\coprod{\rm id}_{F(C)}$), and therefore that the morphisms corresponding to all the different paths form $F(A)$ to $F\big((A\square B)\square C\big)$ coincide. Since we can get similar diagrams with $F(B)$ and $F(C)$, it follows from the universal property of the triple coproduct $F(A)\coprod \big(F(B)\coprod F(C)\big)$ that there exists only a unique morphism from $F(A)\coprod \big(F(B)\coprod F(C)\big)$ to $F\big((A\square B)\square C\big)$ and therefore that the whole diagram commutes.

The commutativity of the two diagrams involving the unit elements can be shown similarly.
\end{proof}

Let $\mathcal{C}$ now be a category of algebraic probability spaces and $F$
the functor that maps a pair $(\mathcal{A},\varphi)$ to the algebra $\mathcal{A}$, i.e., that ``forgets'' the linear functional $\varphi$. Suppose that $\mathcal{C}$ is equipped with a tensor product $\square$ with inclusions and that $F(\mathcal{C})$ has a coproduct $\coprod$. Let $(\mathcal{A},\varphi)$, $(\mathcal{B},\psi)$ be two algebraic probability spaces in $\mathcal{C}$, we will denote the pair $(\mathcal{A},\varphi)\square(\mathcal{B},\psi)$ also by $(\mathcal{A}\square\mathcal{B},\varphi\square\psi)$. By Proposition \ref{prop-forget} we have morphisms $F_2(\mathcal{A},\mathcal{B}):\mathcal{A}\coprod\mathcal{B} \to\mathcal{A}\square\mathcal{B}$ that define a natural transformation from the bifunctor $\coprod$ to the bifunctor $\square$. With these morphisms we can define a new tensor product $\widetilde{\square}$ with inclusions by
\[
(\mathcal{A},\varphi)\widetilde{\square}(\mathcal{B},\psi)= \left(\mathcal{A}\coprod \mathcal{B},(\varphi\square\psi)\circ F_2(\mathcal{A},\mathcal{B})\right).
\]
The inclusions are those defined by the coproduct.

\begin{proposition}
If two random variables $f_1:(\mathcal{A}_1,\varphi_1)\to(\mathcal{B},\psi)$ and
$f_1:(\mathcal{A}_1,\varphi_1)\to(\mathcal{B},\psi)$ are independent with respect to $\square$, then they are also independent with respect to $\widetilde{\square}$.
\end{proposition}
\begin{proof}
If $f_1$ and $f_2$ are independent with respect to $\square$, then there exists a random variable $h:(\mathcal{A}_1\square \mathcal{A}_2,\varphi_1\square\varphi_2)\to (\mathcal{B},\psi)$ that makes diagram (\ref{diag-indep2}) in Definition \ref{def-independence} commuting. Then $h\circ F_2(\mathcal{A}_1,\mathcal{A}_2):(\mathcal{A}_1\coprod \mathcal{A}_2,\varphi_1\tilde{\square}\varphi_2)\to (\mathcal{B},\psi)$ makes the corresponding diagram for $\tilde{\square}$ commuting.
\end{proof}
The converse is not true. Consider the category of algebraic probability spaces with the tensor product, see Subsection \ref{tensor indep}, and take $B=\mathcal{A}_1\coprod\mathcal{A}_2$ and $\psi=(\varphi_1\otimes\varphi_2)\circ F_2(\mathcal{A}_1,\mathcal{A}_2)$. The canonical inclusions $i_{\mathcal{A}_1}:(\mathcal{A}_1,\varphi_1)\to (\mathcal{B},\psi)$ and $i_{\mathcal{A}_2}:(\mathcal{A}_2,\varphi_2)\to (\mathcal{B},\psi)$ are independent w.r.t.\ $\tilde{\otimes}$, but not with respect to the tensor product itself, because their images do not commute in $\mathcal{B}=\mathcal{A}_1\coprod\mathcal{A}_2$.

We will call a tensor product with inclusions in a category of quantum probability spaces {\em semi-universal}, if it is equal to the coproduct of the corresponding category of algebras on the algebras. The preceding discussion shows that every tensor product on the category of algebraic quantum probability spaces $\eufrak{AlgProb}$ has a quasi-universal version.

\section{The Classification of Independences in the Category of Algebraic Probability Spaces}\label{class-indep}

We will now reformulate the classification by Muraki\cite{muraki02} and by Ben Ghorbal and Sch\"urmann\cite{benghorbal01,benghorbal+schuermann99} in terms of semi-universal tensor products with inclusions for the category of algebraic probability spaces $\eufrak{AlgProb}$.

In order to define a semi-universal tensor product with inclusions on $\eufrak{AlgProb}$ one needs a map that associates to a pair of unital functionals $(\varphi_1,\varphi_2)$ on two algebras $\mathcal{A}_1$ and $\mathcal{A}_2$ a unital functional $\varphi_1\cdot\varphi_2$ on the free product $\mathcal{A}_1\coprod\mathcal{A}_2$ (with identification of the units) of $\mathcal{A}_1$ and $\mathcal{A}_2$ in such a way that the bifunctor
\[
\square:(\mathcal{A}_1,\varphi_1)\times(\mathcal{A}_2,\varphi_1)\mapsto (\mathcal{A}_1\coprod\mathcal{A}_2,\varphi_1\cdot\varphi_2)
\]
satisfies all the necessary axioms. Since $\square$ is equal to the coproduct $\coprod$ on the algebras, we don't have a choice for the isomorphisms $\alpha,\lambda,\rho$ implementing the associativity and the left and right unit property, we have to take the ones following from the universal property of the coproduct. The inclusions and the action of $\square$ on the morphisms also have to be the ones given by the coproduct.

The associativity gives us the condition
\begin{equation}\label{cond-associativity}
\big((\varphi_1\cdot\varphi_2)\cdot\varphi_3\big)\circ \alpha_{\mathcal{A}_1,\mathcal{A}_2,\mathcal{A}_3} = \varphi_1\cdot(\varphi_2\cdot\varphi_3),
\end{equation}
for all $(\mathcal{A}_1,\varphi_1),(\mathcal{A}_2,\varphi_2),(\mathcal{A}_3,\varphi_3)$ in $\eufrak{AlgProb}$. Denote the unique unital functional on $\mathbb{C}\mathbf{1}$ by $\delta$, then the unit properties are equivalent to
\[
(\varphi\cdot \delta)\circ\rho_\mathcal{A} = \varphi\quad\mbox{ and }\quad (\delta\cdot\varphi)\circ\lambda_\mathcal{A} = \varphi,
\]
for all $(\mathcal{A},\varphi)$ in $\eufrak{AlgProb}$.
The inclusions are random variables, if and only if
\begin{equation}\label{cond-inclusion}
(\varphi_1\cdot\varphi_2)\circ i_{\mathcal{A}_1}=\varphi_1\quad\mbox{ and }\quad(\varphi_1\cdot\varphi_2)\circ i_{\mathcal{A}_2}=\varphi_2
\end{equation}
for all $(\mathcal{A}_1,\varphi_1),(\mathcal{A}_2,\varphi_2)$ in $\eufrak{AlgProb}$. Finally, from the functoriality of $\square$ we get the condition
\begin{equation}\label{cond-functoriality}
(\varphi_1\cdot\varphi_2)\circ (j_1\coprod j_2) = (\varphi_1\circ j_1)\cdot(\varphi_2\circ j_2)
\end{equation}
for all pairs of morphisms $j_1:(\mathcal{B}_1,\psi_1)\to(\mathcal{A}_1,\varphi_1)$, $j_2:(\mathcal{B}_2,\psi_2)\to(\mathcal{A}_2,\varphi_2)$ in $\eufrak{AlgProb}$.

Our Conditions (\ref{cond-associativity}), (\ref{cond-inclusion}), and (\ref{cond-functoriality}) are exactly the axioms (P2), (P3), and (P4) in Ben Ghorbal and Sch\"urmann\cite{benghorbal+schuermann99}, or the axioms (U2), the first part of (U4), and (U3) in Muraki\cite{muraki02}.
\begin{theorem}
(Muraki\cite{muraki02}, Ben Ghorbal and Sch\"urmann\cite{benghorbal01,benghorbal+schuermann99}) There exist exactly two semi-universal tensor products with inclusions on the category of algebraic probability spaces $\eufrak{AlgProb}$, namely the semi-universal version $\tilde{\otimes}$ of the tensor product defined in Section \ref{tensor indep} and the one associated to the free product $*$ of states.
\end{theorem}

Voiculescu's\cite{voiculescu+dykema+nica92} free product $\varphi_1*\varphi_2$ of two unital functionals can be defined recursively by
\[
(\varphi_1*\varphi_2)(a_1a_2\cdots a_m)= \sum_{I\subsetneqq \{1,\ldots,m\}} (-1)^{m-\sharp I+1} (\varphi_1*\varphi_2)\left(\prod_{k\in I}^{\rightarrow}a_k\right) \prod_{k\not\in I} \varphi_{\epsilon_k}(a_k)
\]
for a typical element $a_1a_2\cdots a_m\in\mathcal{A}_1\coprod\mathcal{A}_2$, with $a_k\in\mathcal{A}_{\epsilon_k}$, $\epsilon_1\not=\epsilon_2\not=\cdots\not=\epsilon_m$, i.e.\ neighboring $a$'s don't belong to the same algebra. $\sharp I$ denotes the number of elements of $I$ and $\prod_{k\in I}^{\rightarrow}a_k$ means that the $a$'s are to be multiplied in the same order in which they appear on the left-hand-side. We use the convention $(\varphi_1*\varphi_2)\left(\prod_{k\in \emptyset}^{\rightarrow}a_k\right)=1$.

Ben Ghorbal and Sch\"urmann\cite{benghorbal01,benghorbal+schuermann99} and Muraki\cite{muraki02} have also considered the category of non-unital algebraic probability $\eufrak{nuAlgProb}$ consisting of pairs $(\mathcal{A},\varphi)$ of a not necessarily unital algebra $\mathcal{A}$ and a linear functional $\varphi$. The morphisms in this category are algebra homomorphisms that leave the functional invariant. On this category we can define three more tensor products with inclusions corresponding to the boolean product $\diamond$, the monotone product $\triangleright$ and the anti-monotone product $\triangleleft$ of states. They can be defined by
\begin{eqnarray*}
\varphi_1\diamond\varphi_2(a_1a_2\cdots a_m) &=& \prod_{k=1}^m \varphi_{\epsilon_k}(a_k), \\
\varphi_1\triangleright\varphi_2(a_1a_2\cdots a_m) &=& \varphi_1\left(\prod_{k:\epsilon_k=1}^{\rightarrow} a_k\right)\prod_{k:\epsilon_k=2} \varphi_2(a_k), \\
\varphi_1\triangleleft\varphi_2(a_1a_2\cdots a_m) &=& \prod_{k:\epsilon_k=1} \varphi_1(a_k)\,\,\varphi_2\left(\prod_{k:\epsilon_k=2}^{\rightarrow} a_k\right), 
\end{eqnarray*}
for $a_1a_2\cdots a_m\in\mathcal{A}_1\coprod\mathcal{A}_2$, $a_k\in\mathcal{A}_{\epsilon_k}$, $\epsilon_1\not=\epsilon_2\not=\cdots\not=\epsilon_m$, i.e.\ neighboring $a$'s don't belong to the same algebra. Note that $\coprod$ denotes here the free product without units, the coproduct in the category of not necessarily unital algebras.

For the classification in the non-unital case, Muraki imposes the additional condition
\begin{equation}\label{cond-factorisation}
(\varphi_1\cdot\varphi_2)(a_1a_2) = \varphi_{\epsilon_1}(a_1)\varphi_{\epsilon_2}(a_2)
\end{equation}
for all $(\epsilon_1,\epsilon_2)\in\big\{(1,2),(2,1)\big\}$, $a_1\in\mathcal{A}_{\epsilon_1}$, $a_2\in\mathcal{A}_{\epsilon_2}$.

\begin{theorem}
(Muraki\cite{muraki02})
There exist exactly five semi-universal tensor products with inclusions satisfying (\ref{cond-factorisation}) on the category of non-unital algebraic probability spaces $\eufrak{nuAlgProb}$, namely the semi-universal version $\tilde{\otimes}$ of the tensor product defined in Section \ref{tensor indep} and the ones associated to the free product $*$, the boolean product $\diamond$, the monotone product $\triangleright$ and the anti-monotone product $\triangleleft$.
\end{theorem}
The monotone and the anti-monotone are not symmetric, i.e.\ $(\mathcal{A}_1\coprod\mathcal{A}_2,\varphi_1\triangleright\varphi_2)$ and $(\mathcal{A}_2\coprod\mathcal{A}_2,\varphi_2\triangleright\varphi_1)$ are not isomorphic in general. Actually, the anti-monotone product is simply the mirror image of the monotone product,
\[
(\mathcal{A}_1\coprod\mathcal{A}_2,\varphi_1\triangleright\varphi_2)\cong(\mathcal{A}_2\coprod\mathcal{A}_1,\varphi_2\triangleleft\varphi_1)
\]
for all $(\mathcal{A}_1,\varphi_1),(\mathcal{A}_2,\varphi_2)$ in the category of non-unital algebraic probability spaces. The other three products are symmetric.

At least in the symmetric setting of Ben Ghorbal and Sch\"urmann, Condition (\ref{cond-factorisation}) is not essential. If one drops it and adds symmetry, one finds in addition the degenerate product
\[
(\varphi_1\cdot \varphi_2) (a_1a_2\cdots a_m) =
\left\{
\begin{array}{lcl}
\varphi_{\epsilon_1}(a_1) & \mbox{ if } & m=1, \\
0 & \mbox{ if } & m>1.
\end{array}
\right.
\]
and a whole family
\[
\varphi_1\bullet_q \varphi_2 = q\big((q^{-1}\varphi_1)\bullet(q^{-1}\varphi_2)\big),
\]
parametrized by a complex number $q\in\mathbb{C}\backslash\{0\}$, for each of the three symmetric products, $\bullet\in\{\tilde{\otimes},*,\diamond\}$. 

\section{The Reduction of Boolean, Monotone, and Anti-Monotone Independence to Tensor Independence}\label{section-reduction-boolean-etc}

We will now present the unification of tensor, monotone, anti-monotone, and boolean independence of Franz\cite{franz02} in our category theoretical framework. It resembles closely the bosonization of Fermi independence in Subsection \ref{reduction-fermi-bose}, but the group $\mathbb{Z}_2$ has to be replaced by the semigroup $M=\{\mathbf{1},p\}$ with two elements, $\mathbf{1}\cdot \mathbf{1}=\mathbf{1}$, $\mathbf{1}\cdot p=p\cdot \mathbf{1}=p\cdot p=p$. We will need the linear functional $\varepsilon:\mathbb{C}M\to\mathbb{C}$ with $\varepsilon(\mathbf{1})=\varepsilon(p)=1$.

The underlying functor and the inclusions are the same for the reduction of the boolean, the monotone and the anti-monotone product. They map the algebra $\mathcal{A}$ of $(\mathcal{A},\varphi)$ to the free product $F(\mathcal{A})=\tilde{\mathcal{A}}\coprod\mathbb{C}M$ of the unitization $\tilde{\mathcal{A}}$ of $\mathcal{A}$ and the group algebra $\mathbb{C}M$ of $M$. For the unital functional $F(\varphi)$ we take the boolean product $\tilde{\varphi}\diamond\varepsilon$ of the unital extension $\tilde{\varphi}$ of $\varphi$ with $\varepsilon$. The elements of $F(\mathcal{A})$ can be written as linear combinations of terms of the form
\[
p^\alpha a_1p\cdots p a_mp^\omega
\]
with $m\in\mathbb{N}$, $\alpha,\omega\in\{0,1\}$, $a_1,\ldots.a_m\in\mathcal{A}$, and $F(\varphi)$ acts on them as
\[
F(\varphi)(p^\alpha a_1p\cdots p a_mp^\omega) = \prod_{k=1}^m \varphi(a_k).
\]
The inclusion is simply
\[
J_{\mathcal{A}}: \mathcal{A}\ni a\mapsto a\in F(\mathcal{A}).
\]
The morphism $F_0:F(\mathbb{C}\mathbf{1})=\mathbb{C}M\to\mathbb{C}\mathbf{1}$ is given by the trivial representation of $M$, $F_0(\mathbf{1})=F_0(p)=\mathbf{1}$. 

The only part of the reduction that is different for the three cases are the morphisms
\[
F_2(\mathcal{A}_1,\mathcal{A}_2):\mathcal{A}_1\coprod\mathcal{A}_2\to F(\mathcal{A}_1)\otimes F(\mathcal{A}_2)=(\tilde{\mathcal{A}}\coprod\mathbb{C}M)\otimes(\tilde{\mathcal{A}}\coprod\mathbb{C}M).
\]
We set
\[
F^{\rm B}_2(\mathcal{A}_1,\mathcal{A}_2)(a) =
\left\{
\begin{array}{lcl}
a\otimes p & \mbox{ if } a\in\mathcal{A}_1, \\
p\otimes a & \mbox{ if } a\in\mathcal{A}_2,
\end{array}
\right.
\]
for the boolean case,
\[
F^{\rm M}_2(\mathcal{A}_1,\mathcal{A}_2)(a) =
\left\{
\begin{array}{lcl}
a\otimes p & \mbox{ if } a\in\mathcal{A}_1, \\
\mathbf{1}\otimes a & \mbox{ if } a\in\mathcal{A}_2,
\end{array}
\right.
\]
for the monotone case, and
\[
F^{\rm AM}_2(\mathcal{A}_1,\mathcal{A}_2)(a) =
\left\{
\begin{array}{lcl}
a\otimes \mathbf{1} & \mbox{ if } a\in\mathcal{A}_1, \\
p\otimes a & \mbox{ if } a\in\mathcal{A}_2,
\end{array}
\right.
\]
for the anti-monotone case.

For the higher order inclusions $J^\bullet_n= F^\bullet_n(\mathcal{A}_1,\ldots,\mathcal{A}_n)\circ J_{\mathcal{A}_1\coprod\cdots\coprod\mathcal{A}_n}$, $\bullet\in\{{\rm B},{\rm M},{\rm AM}\}$, one gets
\begin{eqnarray*}
J_n^{\rm B}(a) &=& p^{\otimes (k-1)}\otimes a \otimes p^{\otimes (n-k)}, \\
J_n^{\rm M}(a) &=& \mathbf{1}^{\otimes (k-1)}\otimes a \otimes p^{\otimes (n-k)}, \\
J_n^{\rm AM}(a) &=& p^{\otimes (k-1)}\otimes a \otimes \mathbf{1}^{\otimes (n-k)},
\end{eqnarray*}
if $a\in\mathcal{A}_k$.

One can verify that this indeed defines reductions $(F^{\rm B},J)$,  $(F^{\rm M},J)$, and $(F^{\rm AM},J)$ from the categories $(\eufrak{nuAlgProb},\diamond,i)$, $(\eufrak{nuAlgProb},\triangleright,i)$, and $(\eufrak{nuAlgProb},\triangleleft,i)$ to $(\eufrak{AlgProb},\otimes,i)$. The functor $U:\eufrak{nuAlgProb}\to\eufrak{AlgProb}$ mentioned at the end of Section \ref{reduction} is the unitization of the algebra and the unital extension of the functional and the morphisms.

This reduces all calculations involving the boolean, monotone or anti-monotone product to the tensor product. These constructions can also be applied to reduce the quantum stochastic calculus on the boolean, monotone, and anti-monotone Fock space to the boson Fock space. Furthermore, they allow to reduce the theories of boolean, monotone, and anti-monotone L\'evy processes to Sch\"urmann's\cite{schuermann93} theory of L\'evy processes on involutive bialgebras, see Franz\cite{franz02}.

\section{Conclusion}

We have seen that the notion of independence in classical and in quantum probability depends on a product structure which is weaker than a universal product and stronger than a tensor product. We gave an abstract definition of this kind of product, which we named tensor product with projections or inclusions, and defined the notion of reduction between these products. We showed how the bosonization of Fermi independence and the reduction of the boolean, monotone, and anti-monotone independence to tensor independence fit into this framework.

We also recalled the classifications of independence by Ben Ghorbal and Sch\"urmann\cite{benghorbal01,benghorbal+schuermann99} and Muraki\cite{muraki02} and showed that their results classify in a sense all tensor products with inclusions on the categories of algebraic probability spaces and non-unital algebraic probability spaces, or at least their semi-universal versions.

There are two ways to get more than the five universal independences. Either one can consider categories of algebraic probability spaces with additional structure, like for Fermi independence, cf.\ Subsection \ref{fermi-indep}, and braided independence, cf. Franz, Schott, and Sch\"urmann\cite{franz+schott+schuermann98}, or one can weaken the assumptions, drop, e.g., associativity, see M{\l}otkowski\cite{mlotkowski99} and the references therein. Romuald Lenczewski\cite{lenczewski98} has given a tensor construction for a family of new products called $m$-free that are not associative, see also Franz and Lenczewski\cite{franz+lenczewski99}. His construction is particularly interesting, because in the limit $m\to\infty$ it approximates the free product. But it is not known, if a reduction of the free product to the tensor product in the sense of Definition \ref{def-reduction} exists.

\end{document}